\documentclass[a4paper]{article}
\usepackage{graphicx}
\usepackage{amsmath}
\usepackage{amsfonts}
\usepackage{amsthm}
\usepackage{amssymb}
\usepackage{subfigure}
\usepackage{multirow}
\usepackage{rotating}
\usepackage{fancybox}
\usepackage{boxedminipage}
\usepackage{version}
\usepackage[colorlinks=true]{hyperref}
\newcommand{\be}{\begin{equation}}
\newcommand{\ee}{\end{equation}}
\newcommand{\ds}{\displaystyle}
\newcommand{\NN}{\mathbb{N}}
\newcommand{\RR}{\mathbb{R}}

\newcommand{\FF}{\mathbb{F}}

\newcommand{\PP}{\mathbb{P}}

\newcommand{\QQ}{\mathbb{Q}}
\newcommand{\A}{\mathbb{A}}

\theoremstyle{definition}
\newtheorem{thm}{Theorem}

\newtheorem{rmk}{Remark}
\newtheorem{dfn}{Definition}

\newtheorem{cor}{Corollary}

\begin{document}
\title{On nowhere continuous Costas functions and infinite Golomb rulers}
\author{Konstantinos Drakakis}
\maketitle

\begin{abstract}
We prove the existence of nowhere continuous bijections that satisfy the Costas property, as well as (countably and uncountably) infinite Golomb rulers. We define and prove the existence of real and rational Costas clouds, namely nowhere continuous Costas injections whose graphs are everywhere dense in a region of the real plane, based on nonlinear solutions of Cauchy's functional equation. We also give 2 constructive examples of a nowhere continuous function, that satisfies a constrained form of the Costas property (over rational or algebraic displacements only, that is), based on the indicator function of a dense subset of the reals. 
\end{abstract}

\section{Introduction}

Costas arrays are square arrangements of dots and blanks (permutation arrays) such that there be exactly one dot per row and column, and such that no 4 dots form a parallelogram and no 3 dots on the same straight line are equidistant. They arose in the 1960s in connection with the development of SONAR/RADAR frequency-hopped waveforms with ideal auto-correlation properties \cite{C1,C2}, but have been the subject of increasingly intensive mathematical study ever since Prof.\ S.\ Golomb published in 1984 \cite{G,GT} some algebraic construction techniques (still the only ones available today) based on finite fields. Mathematicians are mainly concerned with the study of properties of Costas arrays, but also with the settlement of the question of their existence for all orders, which, despite all efforts, still remains open.

Golomb rulers are 1D analogs of Costas arrays: they are linear arrangements of dots and blanks such that no distance between pairs of dots is repeated. The term ``ruler'' arises from the equivalent visualization as a ruler with markings at the integers, where the dots correspond to those markings that get selected. Golomb rulers are, in fact, older than Costas arrays themselves, originally described as Sidon sets, namely sets of integers whose pairwise sums are all distinct (a moment's reflection shows that the 2 definitions are equivalent).

Costas arrays tend to be very irregular: large arrays look like ``clouds'' of dots scattered all over. Similarly, the density of dots in Golomb rulers appears to be very uneven. Naturally, when one contemplates possible generalizations of Costas arrays to the continuum, e.g.\ Costas bijections over an interval of the reals (or the rationals), say $[0,1]$, one's first instinct would most likely dictate that these functions, if they exist at all, will be very irregular. It may then come as a surprise that not only the most obvious such generalizations are actually extremely regular, even infinitely smooth, but also that irregular examples are non-trivial to come by \cite{DR}.

In this work we investigate nowhere continuous functions satisfying the Costas property, and even prove the existence of functions of this type whose graphs are, in addition, everywhere dense in a region of a real plane, taking advantage of discontinuous solutions of Cauchy's integral equation. We also prove the existence of (both countably and uncountably) infinite Golomb rulers which are dense in a (possibly infinite) interval of the real line. Finally, we explicitly construct 2 examples of nowhere continuous Costas bijections, using nowhere continuous indicator functions of dense subsets of the real line as building blocks.

\section{The Costas property}

In this section we provide an overview of the Costas property and, at the same time, an appropriate generalization to non-finite sets.

\subsection{Definition}

Usually the Costas property is defined on the set of the fist $n$ integers $[n]=\{1,2,\ldots,n\}$ \cite{D}. It is necessary to extend this definition in a natural way to cope with (uncountably) infinite sets:
\begin{dfn}\label{cosdef} Let $A,B,T$ be sets with the following properties:
\begin{enumerate}
  \item $A,T$ are subsets of an additive group $G_1$; $B$ is a subset of an additive group $G_2$ ($G_1=G_2$ is allowed).
  \item $0\notin T$.
  \item $\forall t\in T,\ (t+A)\cap A\neq \emptyset$.
\end{enumerate}
A function $f: A\rightarrow B$, so that $f(A)=B$, will be said to have the Costas property from $A$ to $B$ (or on $A$ if $A=B$) with respect to $T$ iff
\[\forall a_1,a_2\in A,\ \forall t\in T: t+a_1,t+a_2\in A,\ \ (f(t+a_1)+f(a_2)=f(t+a_2)+f(a_1)\Rightarrow a_1=a_2)\]
\end{dfn}

\begin{rmk} Historically, the Costas property has been associated exclusively with bijections: this is because the original engineering application that introduced Costas permutations benefits from (but still does not require) bijectivity \cite{C1,C2}. Though, in this work, bijectivity is not considered to be a formally required component of the Costas property, and neither is the weaker condition of injectivity, we will try to respect tradition and focus on Costas functions that are at least injective, and, whenever possible, bijective.
\end{rmk}

\begin{dfn} A Costas permutation of order $n$ is a bijection satisfying Definition \ref{cosdef} with $A=B=[n]$ and $T=[n-1]$ \cite{C1,C2,D}.
\end{dfn}
We map Costas arrays bijectively to Costas permutations using the convention that the $i$th element of the permutation indicates the position of the dot in the $i$th column of the array. Henceforth the terms ``Costas array'' and ``Costas permutation'' will be used interchangeably.

\subsection{Construction of Costas permutations}

Two algebraic construction methods exist for Costas permutations \cite{D,G,GT}. Let us review them briefly without proof, as they will be needed later:

\subsubsection{The Welch construction}

\begin{thm}[Welch construction $W_1(p,\alpha,c)$]\label{w1} Let $p$ be a prime, let $\alpha$ be a primitive root of the finite field $\FF(p)$ of $p$ elements, and let $c\in[p-1]-1$ be a constant; then, the function $f:[p-1]\rightarrow [p-1]$ where $\ds f(i)=\alpha^{i-1+c}\bmod p$ is a bijection with the Costas property.
\end{thm}
The reason for the presence of $-1$ in the exponent is that, when $c=0$, 1 is a fixed point: $f(1)=1$. We refer to arrays generated with $c\neq0$ as circular shifts of the array generated by $c=0$ for the same $p$ and $\alpha$.

\subsubsection{The Golomb construction}

\begin{thm}[Golomb construction $G_2(p,m,a,b)$]\label{g2} Let $p$ be a prime, $m\in\NN$, and let $\alpha,\ \beta$ be primitive roots of the finite field $\FF(p^m)$ of $q=p^m$ elements; then, the function $f:[q-2]\rightarrow [q-2]$ where $\ds \alpha^{i}+\beta^{f(i)}=1$ is a bijection with the Costas property.
\end{thm}

\section{Golomb rulers/ Sidon sets}

Golomb rulers, though bearing the name of Prof.\ S.\ Golomb, were originally described by W.\ C.\ Babcock in the context of an application in telecommunications \cite{B}. It later emerged that they had made their appearance even earlier in the context of harmonic analysis, in a completely equivalent formulation , known as Sidon sets \cite{S}. In this work, the 2 terms will be used interchangeably. A very comprehensive source of information about Golomb rulers is \cite{D4}.

\subsection{Definition}

\begin{dfn}
Let $G$ be an additive group and let $S$ be a subset of $G$ (not necessarily a subgroup). $S$ is a \emph{Sidon set} or a \emph{Golomb ruler} iff
\be \forall x_1,x_2,x_3,x_4\in S:\ x_1+x_2=x_3+x_4\Leftrightarrow \{x_1,x_2\}=\{x_3,x_4\}\ee
\end{dfn}

Clearly, for any $g\in G$, $S$ is a Golomb ruler iff $g+S$ is one. Assuming $G=\NN$, and that $S$ is a set of $m$ integers, we may always assume the first marking to lie at 1; if the last marking lies at $n$, what is the relation between $m$ and $n$? In particular, what is the smallest $n$ for a given $m$? What is the largest $m$ for a given $n$? Golomb rulers that satisfy either of these conditions are dubbed \emph{optimal}. It is conjectured that optimal Golomb rulers asymptotically satisfy the condition
\be\label{grc} m\approx \sqrt{n}.\ee
It is customary to refer to $n$ and $m$ as the \emph{length} and the \emph{number of markings} of the Golomb ruler, respectively.

\subsection{Construction methods for finite Golomb rulers}

Though they won't be needed further in this work, and for the sake of completeness only, we also present some construction methods for Golomb rulers. Note that, contrary to the case of Costas arrays, where the order specifies the number of dots, the definition of a Golomb ruler does not relate the length to the number of markings in any way. Needless to say, construction methods for Golomb rulers considered to be of interest tend to produce reasonably densely populated Golomb rulers, and, in particular, families that asymptotically satisfy (\ref{grc}).

\subsubsection{Erd\"os-Turan construction \cite{ET}}

\begin{thm}
For every prime $p$, the sequence
\be 2pk + (k^2\bmod p),\ k\in[p]-1\ee
forms a Golomb ruler.
\end{thm}

The approximate asymptotic length of such a Golomb ruler with $p$ markings is $2p^2$, hence it deviates from optimality by a factor of 2.

\subsubsection{Rusza-Lindstr\"om construction \cite{L,R2}}

\begin{thm}
Let $p$ be prime, $g$ a primitive root of $\FF(p)$, and $s$ relatively prime to $p-1$. The following sequence
\be (psk + (p-1)g^k) \bmod p(p-1),\ k \in [p]-1\ee
forms a Golomb ruler.
\end{thm}

These rulers are of length (at most) $p(p-1)$ and have $p$ markings, hence they are optimal.

\subsubsection{Bose-Chowla construction \cite{B2,BC}}

\begin{thm}
Let $q = p^n$ be a power of a prime and $g$ a primitive root in $\FF(q^2)$. Then the $q$ integers
\be S=\{i \in [q^2-2]: g^i-g \in \FF(q)\} \ee
have distinct pairwise differences modulo $q^2-1$.

In addition, the set of $q(q-1)$ pairwise differences in $S$, reduced modulo $q^2-1$, equals the set of are all nonzero integers less than $q^2-1$ which are not divisible by $q+1$.
\end{thm}

These rulers are optimal: they have $q$ markings and their length is (at most) $q^2-1$.

\subsubsection{An always applicable construction \cite{D4}}

The previous constructions work only when the number of markings is a (power of a) prime. The following construction works always, but, unfortunately, is far from optimal.

\begin{thm}
For any $n\in\NN^*$, and for a fixed $a\in\{1,2\}$, the sequence
\be ank^2 + k,\  k \in [n]-1\ee
forms a Golomb ruler.
\end{thm}

This ruler has $n$ markings and its length is asymptotically $an^3$, hence it is far from optimal.

\section{Explicit constructions of nowhere continuous Costas functions}
The following two examples use the indicator function of a dense subset $S$ of $\RR$ as a building block, denoted by $\mathbf{1}_S$; specifically, $\mathbf{1}_S(x)=1$ if $\in S$ and $\mathbf{1}_S(x)=0$ otherwise.

\begin{thm}
Let $f:\RR_+\rightarrow \RR_+$ so that
\[f(x)=x^n(1+a\mathbf{1}_\QQ(x)),\ n=2,\ 1+a=c^n,\ c\in\QQ_+.\]
Then $f$ is a nowhere continuous Costas bijection on $\RR_+$ with respect to $\QQ_+$.
\end{thm}

\begin{proof}
We will denote the set of irrational numbers by $\A$. First, we show that $f$ is injective; we distinguish 2 cases for the equation $f(x)=f(y)$:
\begin{itemize}
  \item $x$ and $y$ are both in $\QQ_+$ or both in $\A_+$: we get $x^n=y^n\Leftrightarrow x=y$.
  \item $x\in \QQ_+$, $y\in \A_+$ (the opposite case is obviously similar): we get $(1+a)x^n=(cx)^n=y^n\Leftrightarrow y=cx$ which is impossible, as this equation implies that $y\in\QQ_+$ as well.
\end{itemize}
Therefore, $f$ is injective.

We now show that $f$ is surjective; we distinguish 2 cases for the equation $f(x)=y\in \RR_+$:
\begin{itemize}
  \item $x\in \QQ_+$: then $(cx)^n=y\Leftrightarrow x=\sqrt[n]{y}/c$.
  \item $x\in \A_+$: then $x^n=y\Leftrightarrow x=\sqrt[n]{y}$.
\end{itemize}
For every $y$, then, each of the 2 cases formally yields a solution, but the injectivity of $f$ guarantees that only one will be admissible: specifically, if $y$ is the $n$th power of a rational, the solution comes from the first case, otherwise from the second. In particular, we find an $x$ for every $y$, therefore $f$ is surjective.

Since $\QQ_+$ is dense in $\RR_+$, the indicator function $\mathbf{1}_\QQ$ is nowhere continuous in $\RR_+$, and therefore $f$ is nowhere continuous as well.

We finally need to show that $f$ has the Costas property. Note that, so far, the assumption that $n=2$ was not needed; it will be needed here. We consider the equation
\be f(x+z)-f(x)=f(y+z)-f(y),\ x,y\in\RR_+,\ z\in\QQ_+^*\ee
and observe that $x+z$ is rational iff $x$ is rational. We only need to consider then 2 cases:
\begin{itemize}
  \item $x$ and $y$ are both in $\QQ_+$ or both in $\A_+$: we get
  \be (x+z)^n-x^n=(y+z)^n-y^n\Leftrightarrow (x-y)\sum_{i=1}^{n-1}z^{n-i}{n\choose i}\sum_{j=0}^{i-1} x^{i-1-j}y^j\Leftrightarrow x=y,\ee
  since the double sum is strictly positive (unless $x=y=0$).
  \item $x\in \QQ_+$, $y\in \A_+$ (the opposite case is obviously similar): we get
  \be c^n[(x+z)^n-x^n]=(y+z)^n-y^n.\ee
  For a fixed $x$ and $z\neq 0$, this is a polynomial in $y$ of degree $n-1$ of rational coefficients, and there is, in general, no reason why it could not have irrational solutions. If, however, $n=2$, the equation reduces to
  \be (y+z)^2-y^2=c^2[(x+z)^2-x^2]\Leftrightarrow y=\frac{c^2(2xz+z^2)-z^2}{2z},\ee
  implying that $y\in\QQ$, a contradiction.
\end{itemize}
Overall then, $f(x+z)-f(x)=f(y+z)-f(y)\Rightarrow x=y$ and the Costas property is confirmed. This completes the proof.
\end{proof}

The previous idea can be successfully generalized to yield an expanded family of nowhere continuous functions with the Costas property if $\QQ$ is substituted by $\PP$, the set of all algebraic numbers:

\begin{dfn}
The set $\PP$ of algebraic real numbers consists of those real numbers that can be solutions of polynomials with rational (or, equivalently, integer) coefficients.
\end{dfn}
Note that this set includes not only the rationals but also many irrational numbers.

\begin{thm}
Let $f:R_+^*\rightarrow R_+^*$ so that
\[f(x)=x^s(1+a\mathbf{1}_\PP(x)),\ s\in\QQ\backslash\{0,1\},\ a\in\PP.\]
Then $f$ is a nowhere continuous Costas bijection on $\RR_+^*$ with respect to $\PP_+^*$.
\end{thm}

\begin{proof}
The set $\RR\backslash\PP=\PP^c$ is known as the set of transcendental numbers. The proof follows closely the previous proof. Throughout the proof it will be important to keep in mind that $\PP$ is conveniently closed under several operations such as addition, multiplication, inversion, and exponentiation by a rational exponent.

The equation $f(x)=f(y)$ leads to the two possibilities $x^s=y^s$ and $(1+a)x^s=y^s$ according to whether both $x$ and $y$ are of the same type (algebraic or not), or of different types (in which case, without loss of generality, we consider $x\in\PP$, $y\in\PP^c$): the former leads to $x=y$ (for $s\neq 0$) while the latter implies $\ds y=x(1+a)^{1/s}$, whereby $y$ is algebraic, a contradiction. It follows that $f$ is injective.

The equation $f(x)=y$ leads to the two possibilities $x^s(1+a)=y$ or $x^s=y$, according to whether $x\in\PP$ or not, respectively: the former yields an admissible solution iff $y$ is algebraic, and the latter iff $y$ is transcendental. Hence, a valid $x$ corresponds to each $y$ and $f$ is surjective.

Since $\PP_+$ is dense in $\RR_+$, the indicator function $\mathbf{1}_\PP$ is nowhere continuous in $\RR_+$, and therefore $f$ is nowhere continuous as well.

Finally, we need to show that $f$ satisfies the Costas property. To see this, we form the equation
\be \label{cth2} f(x+z)-f(x)=f(y+z)-f(y),\ x,y\in\RR_+^*,\ z\in\PP_+^*\ee
and consider the resulting 2 cases:
\begin{itemize}
  \item If $x$ and $y$ are both in $\PP_+$ or both in $(\PP^c)_+$, (\ref{cth2}) becomes
  \be (x+z)^s-x^s=(y+z)^s-y^s.\ee
  According to Theorem 5 in \cite{DR}, since the function $g(x)=x^s$ is strictly monotonic with a strictly monotonic derivative in $\RR_+$, it has the Costas property over this set, hence the only possible solution is $x=y$.
  \item If $x$ and $y$ are of different types, let (without loss of generality) $x\in\PP_+$ and $y\in(\PP^c)_+$: the equation becomes
  \be (1+a)[(x+z)^s-x^s]=(y+z)^s-y^s.\ee
  For a fixed $x$, $z\neq 0$, and $s\neq 1$, this is an algebraic equation in $y$, implying that $y$ is algebraic, a contradiction.
\end{itemize}
Overall then, $f(x+z)-f(x)=f(y+z)-f(y)\Rightarrow x=y$ and the Costas property is confirmed. This completes the proof.
\end{proof}
Note that $\RR^*_+$ was chosen instead of $\RR_+$ as the domain of $f$ in e preceding theorem because $f$ is not defined at 0 when $s<0$; $\RR_+$ could have been used for $s>0$.

Although both constructions above are nowhere continuous, their behavior is not as ``wild'' as one might have hoped for: their graph is entirely included within 2 smooth curves, given by the equations $y=x^s$ and $y=Kx^s$, $K>1$, the latter containing countably many points. Is there a function satisfying the Costas property whose graph is everywhere dense on a region of the real plane? This leads to the notion of Costas clouds, studied below.

\section{Costas clouds}

The Costas cloud is a nowhere continuous function (defined on an appropriate set), with an everywhere dense graph, that satisfies the Costas property. More specifically, we propose the following
\begin{dfn}\label{cldef}
Let $f: A\rightarrow B$ be an injection, $A$ being a finite or infinite interval of $\RR$ (such as $[0,1]$, $\RR_+$, $\RR_+^*$, $\RR$ etc.) $f$ will be called a \emph{Costas cloud} on $A\times B$ iff
\begin{enumerate}
  \item It is nowhere continuous on $A$, that is
  \be \forall x\in A, \exists \epsilon>0:\ \forall \delta>0, (\forall y\in A: |x-y|<\delta)\rightarrow  |f(x)-f(y)|>\epsilon.\ee
  \item The graph of $f$ is everywhere dense in $A\times B$, that is
  \be \forall x,y\in A,\forall \epsilon>0,\ \ \exists z\in A:\ |z-x|<\epsilon,\ |f(z)-y|<\epsilon.\ee
  \item $f$ has the Costas property from $A$ to $B$ with respect to $T=\{t\in\RR^*: (t+A)\cap A\neq \emptyset\}$.
\end{enumerate}
If, in addition, $f$ is bijective, it will be a \emph{bijective Costas cloud}.
\end{dfn}
Note that $t\in T\Leftrightarrow -t\in T$, so we can confine our attention to $t>0$. Our intention is to prove that bijective Costas clouds exist; to carry out the proof, we will need some background.

\subsection{Cauchy's functional equation}

\begin{dfn} Let $f:\RR\rightarrow \RR$; it satisfies \emph{Cauchy's functional equation} iff
\be \forall x,y\in\RR, f(x+y)=f(x)+f(y). \label{ce}\ee
\end{dfn}
A detailed study of this equation can be found in \cite{AD}. For the sake of completeness, we state and prove below those properties of the solutions we will need.

\begin{thm}\label{cepr}
The solution $f$ of Cauchy's equation (\ref{ce}) satisfies the following properties:
\begin{enumerate}
  \item $\forall q\in\QQ,\forall x\in\RR\ f(qx)=qf(x)$.
  \item $f$ is continuous iff $\exists c\in\RR:\ \forall x\in \RR\ f(cx)=cx$.
  \item $f$ is continuous everywhere iff it is continuous at a point.
  \item $f$ is discontinuous iff its graph is everywhere dense on the real plane.
\end{enumerate}
\end{thm}

\begin{proof}\
\begin{itemize}
  \item Setting $x=y=0$ we get $f(0+0)=f(0)=f(0)+f(0)\Leftrightarrow f(0)=0$.
  \item Setting $y=-x$ we get $f(x-x)=f(0)=0=f(x)+f(-x)\Leftrightarrow f(-x)=-f(x)$.
  \item Setting $x=x_1$, $y=x_2+\ldots+x_n$ for $n\in\NN^*$ we get $f(x_1+x_2+\ldots+x_n)=f(x_1)+f(x_2+\ldots+x_n)=\ldots=f(x_1)+f(x_2)+\ldots+f(x_n)$.
  \item Setting $x_1=x_2=\ldots=x_n=x$ we get $f(nx)=nf(x)$. Setting $y=nx$ we get $\ds f\left(\frac{y}{n}\right)=\frac{1}{n}f(y)$.
\end{itemize}
Expressing the rational $q$ as $\ds \frac{m}{n}$, all of the above shows that, for any $x\in\RR$:
\be f(qx)=f\left(\frac{m}{n}x\right)=(-1)^{\text{sign}(m)}f\left(\frac{|m|}{n}x\right)=(-1)^{\text{sign}(m)}|m|f\left(\frac{1}{n}x\right)=\frac{m}{n}f(x)=qf(x).\ee

Assume now $f$ is continuous: for every $x\in \RR$ there exists a sequence $\{q_n\}$ of rationals such that $q_n\rightarrow x$. It follows that
\be f(x)=f(\lim q_n)=\lim f(q_n)=\lim q_n f(1)=xf(1)=cf(x),\ c=f(1).\ee
Conversely, every function of the form $f(x)=cx$ satisfies Cauchy's equation (\ref{ce}) and is continuous.

Assume, without loss of generality, that $f$ is continuous at 0; then
\be \lim_{y\rightarrow 0} f(x+y)=f(x)+\lim_{y\rightarrow 0} f(y)=f(x)+f(0)=f(x),\ee
hence $f$ is continuous at (an arbitrary) $x$.

Assume now $f$ is not continuous: then, by what we just proved, it must be nowhere continuous and it cannot be linear (though it has to be linear over the rationals). This implies
\be \exists x_1,x_2\in \RR: \frac{f(x_1)}{x_1}\neq \frac{f(x_2)}{x_2},\ee
whence it follows that the 2 vectors $v_1=(x_1,f(x_1))$, $v_2=(x_2,f(x_2))$ are linearly independent and, consequently, span the entire real plane. The set of vectors, then, of the form $\{r_1v_1+r_2v_2: r_1,r_2\in\QQ\}$ are an everywhere dense subset of the real plane; but
\be r_1v_1+r_2v_2=r_1(x_1,f(x_1))+r_2(x_2,f(x_2))=(r_1x_1+r_2x_2,f(r_1x_1+r_2x_2)),\ee
which means that the subset $\{f(x): x=r_1x_1+r_2x_2,\ r_1,r_2\in\QQ\}$ of the graph of $f$ is everywhere dense on the plane, whence the graph of $f$ itself is everywhere dense on the plane. Conversely, if the graph of $f$ is everywhere dense on the real plane, $f$ cannot possibly be continuous, or else it would be linear and thus would not possess an everywhere dense graph. This completes the proof.
\end{proof}

\begin{rmk}
Cauchy's functional equation, despite its simplicity, has been playing a prominent role in analysis: Hilbert's 5th problem essentially proposes a generalization of this equation, while an important area of study is the Hyers-Rassias-Ulam stability of this equation (or slight variants thereof) \cite{H,HR,R}.
\end{rmk}

\subsection{Bijective solutions of Cauchy's functional equation}

\begin{thm}\label{cesnc}
There exist solutions of Cauchy's equation that are nowhere continuous bijections, everywhere dense on the real plane.
\end{thm}

\begin{proof}
Consider $\RR$ as a vector space over $\QQ$. This vector space must necessarily have an uncountable basis, or else $\RR$ itself would be countable: it follows by the Continuum Hypothesis that this basis can be indexed by the real numbers, and, therefore, that we can describe the basis as $B=\{b_a: a\in\RR\}$. By definition, any real number $x$ admits a (finite) linear expansion $x=q_1b_{a_1}+\ldots+q_nb_{a_n}$ over this basis (where the rational $q$, the indices $a$, and $n$ are obviously functions of $x$). Furthermore, a solution $f$ of Cauchy's equation (\ref{ce}) can be considered as a linear map over this vector space, and we can write
\be f(x)=f(q_1b_{a_1}+\ldots+q_nb_{a_n})=q_1f(b_{a_1})+\ldots+q_nf(b_{a_n}),\ee
which verifies the well-known result that a linear map over a vector space is unambiguously defined by its effect on the vector space basis $B$. Assuming now that $\exists c\in\RR:\ \forall a\in\RR,\ f(b_a)=cb_a$, we get $f(x)=cx$ for all $x$, namely that $f$ is linear. By Theorem \ref{cepr}, though, this will be the \emph{only} case resulting to a linear $f$: in all other cases $f$ will be nowhere continuous and its graph everywhere dense. Choosing a bijection $g$ over $\RR$ (other than the identity) such that $\forall a\in \RR, f(b_a)=b_{g(a)}$ results to an $f$ that is bijective as well. This completes the proof.
\end{proof}

\begin{rmk}\label{qp}
The theorem does not rely on the exact nature of $\QQ$: it only requires a field over which $\RR$ has an uncountable basis. In particular, any countable field extension of $\QQ$ would have been equally suitable, such as $\PP$.
\end{rmk}

\subsection{The existence of bijective Costas clouds on $\RR\times \RR^*_+$ --- the Welch method}

\begin{thm}\label{clthm}
Bijective Costas clouds on $\RR\times \RR^*_+$ exist.
\end{thm}

\begin{proof} Consider a nowhere continuous bijection $f$ on $\RR$ whose graph is everywhere dense on $\RR^2$, and that further satisfies Cauchy's equation. Construct $g: \RR\rightarrow \RR_+^*$ such that $\ds g(x)=\exp(f(x))$; $g$ inherits the properties $f$ has, and is then itself a nowhere continuous bijection (as the exponential function is strictly monotonic), whose graph is everywhere dense in $\RR\times \RR^*_+$. Furthermore, it satisfies the Costas property from $\RR$ to $\RR^*_+$ with respect to $\RR^*_+$:
\begin{multline}
g(x+z)-g(x)=g(y+z)-g(y)\Leftrightarrow \exp(f(x+z))-\exp(f(x))=\exp(f(y+z))-\exp(f(y))\Leftrightarrow\\ (\exp(f(z))-1)(\exp(f(x))-\exp(f(y)))=0\Leftrightarrow \exp(f(z))=1\text{ or } \exp(f(x))=\exp(f(y))\Leftrightarrow\\
f(z)=0\text{ or } f(x)=f(y)\Leftrightarrow z=0\text{ (rejected)}\text{ or } x=y.
\end{multline}
This completes the proof.
\end{proof}

\begin{cor}
Costas clouds on $(\RR^*_+)^2$ exist.
\end{cor}

It has been suggested that the success of the Welch method (Theorem \ref{w1}) in the construction of Costas permutations lies in the interplay between an additive and a multiplicative structure \cite{G2}: indeed, the exponent is an additive function, while the exponentiation turns this additive structure into a multiplicative one. But, according to the proof of Theorem \ref{clthm}, $g(x)=\exp(f(x))$, where $f:\RR\rightarrow \RR$ satisfies Cauchy's functional equation (\ref{ce}) and is additive: this function, therefore, exhibits the same interplay between the additive and the multiplicative structure, and can accordingly be considered a generalization of the Welch construction. Theorem \ref{clthm}, then, viewed from the point of view of the method it follows instead of the result it achieves, reads as follows:

\begin{thm}[Generalized Welch construction in the continuum]
Let $f:\RR\rightarrow \RR$ be a bijection satisfying (\ref{ce}) and let $g:\RR\rightarrow \RR^*_+$ be such that $g(x)=\exp(f(x))$; then $g$ is a bijection satisfying the Costas property with respect to $\RR_+^*$. In particular, if $f$ has an everywhere dense graph in $\RR^2$, so does $g$ in $\RR\times \RR^*_+$.
\end{thm}

\subsection{The existence of bijective Costas clouds on $\RR\times (0,1)$}

Based on Theorem \ref{clthm}, we can apply a simple transformation and obtain a Costas cloud on $\RR\times (0,1)$:

\begin{thm}\label{strclthm}
Bijective Costas clouds on $\RR\times (0,1)$ exist.
\end{thm}

\begin{proof}
Consider the function $g$ of (the proof of) Theorem \ref{clthm}, note that $g(x+y)=g(x)g(y)$, and consider $h:\RR\rightarrow (0,1)$ such that $\ds h(x)=\exp(-g(x))$. It clearly follows that $h$ is bijective, nowhere continuous, and that it has an everywhere dense graph, as $g$ has these properties. We need to verify the Costas property:
\begin{multline} \exp(-g(x+z))-\exp(-g(x))=\exp(-g(y+z))-\exp(-g(y))\Leftrightarrow\\ \exp(-g(x)g(z))-\exp(-g(x))=\exp(-g(y)g(z))-\exp(-g(y))\Leftrightarrow u^a-u=v^a-v,\end{multline}
where $u=\exp(-g(x))$, $v=\exp(-g(y))$, $a=-g(z)$, $u,v\in (0,1)$, $a\in(-\infty,0)$. But $(u^a-u)'=au^{a-1}-1<0$ since $a<0$, whence
\be u^a-u=v^a-v\Leftrightarrow u=v\Leftrightarrow g(x)=g(y)\Leftrightarrow x=y.\ee
This completes the proof.
\end{proof}

\subsection{The existence of Costas clouds on $(0,1)^2$}

As an immediate corollary of Theorem \ref{strclthm} we obtain:

\begin{cor}\label{cor1}
Costas clouds on $(0,1)^2$ exist.
\end{cor}

But we can also apply another simple transformation, again based on Theorem \ref{clthm}, to obtain this result:

\begin{proof}[Alternative proof of Corollary \ref{cor1}]
Consider the function $g$ of (the proof of) Theorem \ref{clthm}, and consider $h:(0,1)\rightarrow (0,1)$ such that $\ds h=\frac{g}{1+g}$. Clearly $h$ is injective, nowhere continuous on $(0,1)$, and has a graph everywhere dense in $(0,1)^2$, because $g$ has all these properties. We only need to make sure it has the Costas property:
\begin{multline}
\frac{g(x+d)}{1+g(x+d)}-\frac{g(x)}{1+g(x)}=\frac{g(y+d)}{1+g(y+d)}-\frac{g(y)}{1+g(y)}\Leftrightarrow
\frac{g(x)g(d)}{1+g(x)g(d)}-\frac{g(x)}{1+g(x)}=\frac{g(y)g(d)}{1+g(y)g(d)}-\frac{g(y)}{1+g(y)}\Leftrightarrow\\
\frac{g(x)(g(d)-1)}{(1+g(x)g(d))(1+g(x))}=\frac{g(y)(g(d)-1)}{(1+g(y)g(d))(1+g(y))}\Leftrightarrow
g(x)(1+g(y)g(d))(1+g(y))=g(y)(1+g(x)g(d))(1+g(x))\Leftrightarrow\\
(g(x)-g(y))(1-g(x)g(y)g(d))=0=(g(x)-g(y))(1-g(x+y+d))\Leftrightarrow x=y\text{ or } x+y+d=0.
\end{multline}
The second alternative is impossible as $x,y,z>0$, so we necessarily obtain $x=y$, and the Costas property is verified. This completes the proof.
\end{proof}

\subsection{Bijective Costas clouds on $(\RR_+^*)^2$ --- the Golomb method}

Since the Welch method can be successfully generalized on the real line, can the same be done for the Golomb method (Theorem \ref{g2})?

\begin{thm}[Generalized Golomb construction in the continuum]
Let $f:\RR\rightarrow \RR$ be a bijection satisfying (\ref{ce}) and let $g:\RR^*\rightarrow \RR^*$ be such that $\ds \exp(g(x))+\exp(f(x))=1$; then $g$ is a bijection satisfying the Costas property with respect to $\RR_+^*$. In particular, if $f$ has an everywhere dense graph in $\RR^2$, so does $g$.
\end{thm}

\begin{proof}
$g$ is clearly bijective iff $f$ is, and it inherits the property of the everywhere dense graph as long as $f$ has it; we just need to show the Costas property:

\begin{multline}
g(x+z)-g(x)=g(y+z)-g(y)\Leftrightarrow\\ \ln(1-\exp(f(x+z)))-\ln(1-\exp(f(x)))=\ln(1-\exp(f(y+z)))-\ln(1-\exp(f(y)))\Leftrightarrow\\
\frac{1-\exp(f(x)f(z))}{1-\exp(f(x))}=\frac{1-\exp(f(y)f(z))}{1-\exp(f(y))}\Leftrightarrow \exp(f(x))\exp(f(z))+\exp(f(y))=\exp(f(y))\exp(f(z))+\exp(f(x))\Leftrightarrow\\ (1-\exp(f(z)))(\exp(f(x))-\exp(f(x)))=0\Leftrightarrow
f(z)=0\text{ or } f(x)=f(y)\Leftrightarrow z=0\text{ (rejected)}\text{ or } x=y.
\end{multline}
Note finally that $x=0\Leftrightarrow f(x)=0\Leftrightarrow \exp(g(0))=0$, which is impossible, so $g$ cannot be defined on $x=0$. This completes the proof.
\end{proof}

Viewed from the point of view of the result instead of the method, the theorem reads:

\begin{thm}
Bijective Costas clouds on $(\RR^*)^2$ exist.
\end{thm}

Though it is really tempting to extend $g$ by setting $g(0)=0$ and fill the ``hole'' on the plane, unfortunately the resulting extension no longer has the Costas property.

\subsection{Rational Costas clouds}

An algorithm to construct Costas bijections on the set $Q=[0,1]\cap \QQ$ has already been proposed in \cite{DR}, but the density of its graph in $[0,1]^2$ was not considered or studied at the time. We propose here a different algorithm that produces a rational Costas cloud. Referring to Definition \ref{cldef}, note that it is of no interest to consider the density of the graph of a function on $A\times A$ whenever $A\subset \QQ$: we correct then this part of the definition to read that the graph has to be dense on $I(A)\times I(A)$, where $I(A)$ denotes the smallest (closed) interval in $\RR$ containing $A$.

The construction proposed in Theorem \ref{clthm} does not produce a rational Costas cloud, as not only the function fails to be a bijection on the rationals, but also the images of all rationals lie on a curve.We need then a new mechanism, such as the one provided by the following

\begin{thm}\label{rcc}
Enumerate $Q$ so that $Q=\{r_n: n\in\NN\}$, and consider two copies of it, $Q_x^1$ and $Q_y^1$ (it is also allowed to use different enumeration schemes in the two sets). Consider the following inductive construction on $S=Q^2$:
\begin{description}
  \item[Stage 1] Draw the horizontal line through $(0,2^{-1})$ and the vertical line through $(2^{-1},0)$, thus dividing $S$ into 4 smaller squares; choose one point with rational coordinates in each so that the Costas property holds, and so that the points $(r_0,u)$ and $(v,r_0)$ are chosen. Let the points be $(x_{1i},y_{1i}),\ i\in[4]$, and set $Q_x^2=Q_x^1-\{x_{1i}:i\in[4]\}$, $Q_y^2=Q_y^1-\{y_{1i}:i\in[4]\}$.
  \item[Stage $n$] Draw the horizontal lines through $(0,k2^{-n}),\ k=1,\ldots,2^n-1$ and the vertical lines through $(k2^{-n},0),\ k=1,\ldots,2^n-1$, thus dividing $S$ into $4^n$ smaller squares. Choose one point $(x_{ni},y_{ni})$ $x_{ni}\in Q_x^n$, $y_{ni}\in Q_y^n$ in each square ($i\in[4^n]$), so that, for all points together chosen in stages 1 through $n$ inclusive, the Costas property holds, and so that, if $N_x(n)$ and $N_y(n)$ are the smallest natural numbers so that $r_{N_x(n)}\in Q^n_x$ and $r_{N_y(n)}\in Q^n_y$, respectively, one of the chosen points has $r_{N_x(n)}$ as its first coordinate and one $r_{N_y(n)}$ as its second. Set $Q_x^{n+1}=Q_x^n-\{x_{ni}:i\in[4^n]\}$, $Q_y^{n+1}=Q_y^n-\{y_{ni}:i\in[4^n]\}$.
\end{description}
The set of points so constructed corresponds to the graph of a bijective rational Costas cloud on $Q^2$.
\end{thm}

\begin{proof} The construction is possible because, for every stage and every square, we are called to choose a rational point in this square so that finitely many constraints are satisfied: this is always possible as there are infinitely many rational points in a square. The resulting function clearly satisfies the Costas property, is injective between its domain and its range, and its graph is everywhere dense in $S$. Further, every rational number in the domain and the range is used exactly once. This completes the proof.
\end{proof}

\begin{rmk} The construction above can be modified to yield bijective rational Costas clouds on the entire $\QQ^2$: just apply stage $n$ to the square $S_n=\{(x,y)\in\RR^2: \max(|x|,|y|)<2^n\}$. In other words, as stages progress, the grid not only becomes more and more refined but also expands.
\end{rmk}

\begin{rmk}
The square grid described in the theorem has the property that the squares in a given stage are all of the same size, and that any 2 squares either share a common boundary or else one contains the other. Neither of these properties is necessary, strictly speaking: for example, we could have used a grid of the form $k17^{-n}$, or even $k\pi^{-n}$.
\end{rmk}

\begin{rmk}
The proof does not depend on the exact nature of $\QQ$, except for the facts that it is countably infinite and everywhere dense in $\RR$: any other set with these properties would have sufficed, such as $\PP$.
\end{rmk}

\section{Countably and uncountably infinite dense Golomb rulers}

Infinite Golomb rulers have been studied in the past \cite{R3}: the main point of interest in these studies has been the behavior of the density $m(n)/n$, $m(n)$ being the number of markings within the set $[n]$. It should be mentioned here that, even though references to ``dense infinite Sidon sequences'' can be found in the literature (for example \cite{AKS}), giving the impression, at first sight, that the results below may not be novel after all, in fact the term ``dense'' is used not in the analytic sense but in conjunction with the behavior of $m(n)/n$. Such references are, then, totally unrelated to the task we are about to embark on.

\subsection{The existence of countably infinite dense Golomb rulers}

\begin{thm}
Let $I$ be a subinterval of $\RR$ (possibly infinite, possibly $\RR$ itself): countably infinite Golomb rulers in $I$ exist, and may be suitably chosen to be dense in $I$.
\end{thm}

\begin{proof}
The idea of the proof is essentially the same as in Theorem \ref{rcc}, and relies upon the enumerability of $\QQ$. Enumerate $Q=\QQ\cap I$ so that $Q=\{r_n: n\in\NN\}$, and then set $S_0=\{r_0\}$. Assume now that $S_m$, containing $m$ rationals, has been created, and that $r_n$ is being currently considered: if $S_m\cup \{r_n\}$ is a Golomb ruler, set $S_{m+1}=S_m\cup \{r_n\}$; proceed to consider $r_{n+1}$. For each $m$ we are bound to construct $S_{m+1}$ out of $S_m$, as there are infinitely many rationals to choose from, and only finitely many constraints due to the Golomb ruler property.  $\ds S=\lim S_m$ is a countably infinite Golomb ruler.

To ensure the ruler is dense in $I$, let $I_{N,k}$ be the $k$th subinterval of length $2^{-N}$ (say counting from left to right) of $I\cap[-2^N,2^N]$, and then order all $I_{N,k}$s consecutively, first by $N$ and then by $k$: $I_{1,1},I_{1,2},I_{1,3},I_{1,4},I_{2,1},\ldots,I_{2,16},\ldots$, thus obtaining a sequence of intervals $I_n$. Now apply the construction proposed above with the extra requirement that, for every $m$, the element $x$ that turns $S_m$ into $S_{m+1}$ must belong in $I_{m+1}$: this is always possible, as $I_{m+1}$ itself contains infinitely many rationals.

This completes the proof.
\end{proof}

\begin{rmk}
Any countably infinite set, such as $\PP$, could have been used in the proof above, instead of $\QQ$; the proof clearly does not rely on the exact nature of the set used.
\end{rmk}

\subsection{The existence of uncountably infinite dense Golomb rulers}

\begin{thm}
Let $I$ be a subinterval of $\RR$ (possibly infinite, possibly $R$ itself): uncountably infinite Golomb rulers in $I$ exist, and may be suitably chosen to be dense in $I$.
\end{thm}

\begin{proof}
Consider $\RR$ to be a vector space $V$ over $\QQ$, and let $B=\{b_a: a\in\RR\}$ be an (uncountable) basis of this vector space (see the proof of Theorem \ref{cesnc}). Consider the family of subspaces $V_x=\text{span}\{b_u:u\leq x\}$, and observe that $V_x\subset V_y$ iff $x\leq y$. Choose a unique point $s_x\in V_x$ such that $s_x\in I$ and that $s_x\notin V_y,\ y<x$, and form the set $S=\{s_x: x\in\RR\}$. This is an uncountable subset of $I$, and we show it is indeed a Golomb ruler. Consider the equation
\be s_{x_1}+s_{x_2}=s_{x_3}+s_{x_4}, x_1<x_2,\ x_3<x_4\Leftrightarrow s_{x_2}=s_{x_3}+s_{x_4}-s_{x_1},\ee
whence $\ds s_{x_2}\in V_{x_2}\cap V_{\max(x_1,x_3,x_4)}$. Since, by assumption, $x_3<x_4$, it follows that $\max(x_1,x_3,x_4)=\max(x_1,x_4)$. If $x_1>x_4$, we obtain $s_{x_2}\in V_{x_1}$ which is impossible; therefore $x_1<x_4$, implying that $x_2=x_4$, and, consequently, that $x_3=x_1$.

Note that, for any $x\in\RR$, any $qs_x\in A$ with $q\in\QQ$ can be chosen instead of $s_x$. This allows us to choose $S$ dense in $I$. This completes the proof.
\end{proof}

\begin{rmk}
$\PP$ could have been used instead of $\QQ$ (see Remark \ref{qp}).
\end{rmk}

\section{Conclusion}

We have constructed explicit examples of nowhere continuous injections satisfying a constrained form of the Costas property (over rational or algebraic displacements only, that is), using the indicator function of the rationals or of the algebraic numbers as a building block. Furthermore, we proved the existence of nowhere continuous injections, based on Cauchy's functional equation, whose graphs are everywhere dense in a region of the real plane; such functions, named Costas clouds, are perhaps what first springs into one's mind when considering possible generalizations of Costas arrays in the continuum, due to the very haphazard positioning of their dots. We considered both real and rational clouds. The real Costas clouds, in particular, led to the generalization in the continuum of the 2 main generation methods for Costas permutations, namely the Welch and the Golomb construction. These functions are highly non-trivial to construct, and their existence non-trivial to prove.

Similarly, we proved the existence of (countably and uncountably) infinite Golomb rulers in a (finite of infinite) interval of the real line, that can optionally be constructed so that they have the extra property of being everywhere dense in this interval. We also noted that, though ``infinite dense'' Golomb rulers have appeared in the literature before, the word ``dense'' had an entirely different meaning and was not used in the analytic sense.

Both uncountable constructions (Costas clouds and Golomb rulers) relied on 2 ideas, namely the consideration of $\RR$ as vector space over $\QQ$ possessing an uncountable basis, and the use of Cauchy's functional equation (in the case of Costas clouds). These ideas not only made the proof of the existence of these objects possible, but also revealed how much freedom we have for their construction; without these ideas the tasks seemed hopeless. Regarding the corresponding countably infinite (e.g.\ rational) objects, enumerability itself is sufficient to both establish existence and allow great freedom of construction.


\section*{Acknowledgements} The author is indebted to Prof.\ Nigel Boston (School of Mathematics, University College Dublin) for his helpful suggestions regarding Theorem \ref{clthm}, as well as Prof.\ Roderick Gow (ibid.) and Dr.\ Scott Rickard (School of Electrical, Electronic \& Mechanical Engineering, University College Dublin) for the many useful discussions on the topic.

\end{document}